\documentclass{article}

\usepackage[latin1]{inputenc}
\usepackage[T1]{fontenc}

\usepackage{latexsym}

\def\R{\relax\ifmmode I\!\!R\else$I\!\!R$\fi}

\def\Z{\relax\ifmmode Z\!\!\!Z\else$Z\!\!\!Z$\fi}

\def\C{\relax\ifmmode C\!\!\!\!I\else$C\!\!\!\!I$\fi}

\def\K{\relax\ifmmode I\!\!K\else$I\!\!K$\fi}

\def\N{\relax\ifmmode I\!\!N\else$I\!\!N$\fi}
\def\H{\relax\ifmmode I\!\!H\else$I\!\!N$\fi}
\def\P{\relax\ifmmode I\!\!P\else$I\!\!N$\fi}

\newtheorem{corollary}{Corollary}[section]

\newtheorem{conjecture}{Conjecture}[section]
\newtheorem{proposition}{Proposition}[section]

\newcommand{\proof}{{\noindent\bf Proof }}

 \newcommand{\qed}{\hfill $\Box$}

\begin{document}
\thispagestyle{empty}
\begin{center}
{\Large {\bf An isomorphism between polynomial eigenfunctions of the transfer operator 
and the Eichler cohomology for modular groups }}
\end{center}
\begin{center}
D.Mayer 
 and J. Neunh\"auserer\\
Institut f\"ur Theoretische Physik,\\ Technische Universit\"at Clausthal.
\footnote{This work is supported by the DFG research group 'Zeta functions 
and locally symmetric spaces'.}\\
Arnold-Sommerfeld-Str. 6, \\ 38678 Clausthal-Zellerfeld, Germany\\
e-mail: dieter.mayer@tu-clausthal.de \\ joerg.neunhaeuserer@tu-clausthal.de
\end{center}
\begin{abstract}
  For the group $PSL(2,\Z)$ it is known that there is an isomorphism
  between polynomial eigenfunctions of the transfer operator for the
  geodesic flow and the Eichler cohomology in the theory of modular forms, see
  \cite{[CM3]}, \cite{[LW]}, \cite{[LZ]}.  In \cite{[CM3]} it is
  indicated that such an isomorphism exists as well for the subgroups
  $\Gamma(2)$ and $\Gamma_{0}(2)$ of $PSL(2,\Z)$.  We will prove this
  and provide some evidence by computer aided algebraic calculations
  that such an isomorphism exists for all principal congruence
  subgroups $\Gamma(N)$ and all congruence subgroups of Hecke type $\Gamma_{0}(N)$.
\end{abstract}
\newpage
\section{Introduction}
The transfer operator $L_{s}$ for the geodesic flow on the surface
corresponding to a modular group is a generalization of the classical
Perron Frobenius operator in ergodic theory for the Poincare map of
this dynamical system. The operator is of special interest because it
provides a new approach to Selberg's zeta function which encodes
information about the length spectrum of the geodesic flow. In fact it
was shown that Selberg's zeta function can be expressed using the
Fredholm determinant of this operator and that the zeros of Selberg's
zeta function correspond to values of $s$ for which the transfer
operator $L_{s}$ has eigenvalue $+1$ or $-1$ (see \cite{[MA]},
\cite{[CM1]}, \cite{[CM2]} and \cite{[CM3]}).  In this paper we are
interested in polynomial eigenfunctions of the transfer operator corresponding to
these eigenvalues. These eigenfunctions are determined by the so called
Lewis equation found in \cite{[CM2]} and \cite{[LZ]}. Our general
conjecture is that using the Lewis equation one can directly find an
isomorphism between the spaces of polynomial eigenfunctions of the
transfer operator and the spaces of Eichler's cohomology classes which
are well known in number theory \cite{[EI]}. In fact these cohomology
classes correspond in a special way to automorphic forms of modular
groups. In the case of $PSL(2,\Z)$ the existence of such an
isomorphism is known and in \cite{[CM3]} we presented numerical
indications that such an isomorphism exists for the simplest subgroups
of $PSL(2,\Z)$. We will present here rigorous arguments proving the
existence of such an isomorphism in the case of the groups $\Gamma(2)$
(see Proposition 7.1) and $\Gamma_{0}(2)$ (see Proposition 8.1 and
Proposition 8.2). Furthermore we will provide strong numerical
evidence that such an isomorphism exists for all principal congruence
subgroups $\Gamma(N)$ (see Table 1) and all congruence subgroups of
Hecke type $\Gamma_{0}(N)$ (see Table 2).  The main problem to prove
this fact rigorously is the complicated combinatorial
interplay between the action of the generators on the coset sets of these groups and the structure of the Lewis equation.   \\
The rest of this paper is organized as follows. In sections 2 and 3 we
give a brief overview about modular groups, modular forms and the
Eichler cohomology emphasizing those facts about congruence subgroups
that we need.  In section 4 we introduce the transfer operator. In
section 5 we state the Lewis equation for modular groups which is then
exploited in section 6 to describe the spaces of polynomial
eigenfunctions of the transfer operator. In section 7 we present our
results for principal congruence modular groups and in section 8 we
present our results for congruence subgroups of Hecke type.

\section{Modular groups and modular surfaces}
We will give here a brief introduction to modular groups and modular
surfaces.  Our references for this material are
the books of Miyake \cite{[MI]} and Shimura \cite{[SCH]}.\\
Consider the upper half plane $\H=\{z\in\C|\Im(z)>0\}$ with the
Poincare metric given by
\[ ds^{2}(z)=\frac{dz~ d\bar z}{\Im(z)^{2}}. \]
The special linear group
\[ SL(2,\R)=\{ \left( \begin{array}{cc}a&b\\ c & d\end{array}  \right)|a,b,c,d\in\R \mbox{ and } ad-bc=1\} \]
acts on the upper half plane $\H$ by the M\"obius transformation
\[\left( \begin{array}{cc}a&b\\ c & d\end{array}  \right) z=\frac{az+b}{cz+d}. \]
The group
\[ PSL(2,\R)=SL(2,\R)/{\pm I}, \]
is known to be the isometry group of $\H$\footnote{$I$ denotes the
  identity matrix through out this paper.}. Discrete subgroups of
$PSL(2,\R)$ like
\[ \Gamma(1)= PSL(2,\Z)=\{ \left( \begin{array}{cc}a&b\\ c & d\end{array}  \right)\in  PSL(2,\R) | a,b,c,d\in\Z\} \]
are called {\bf Fuchsian groups} and the subgroups of $\Gamma(1)$ of
finite index are called {\bf modular groups}.  Given a modular group
$\Gamma$ we denote by $\tilde\Gamma$ the
set $\Gamma\backslash\Gamma(1)$ of right cosets.\\
In this paper we are especially interested in {\bf principal
  congruence subgroups}
\[ \Gamma(N):=\{A\in\Gamma(1)|A = I \mbox{ mod } N\}\]
and {\bf congruence subgroups of Hecke type}
\[ \Gamma_{0}(N):= \{ \left( \begin{array}{cc}a&b\\ c & d\end{array}  \right)\in  \Gamma(1)| c=0 \mbox{ mod } N \} \]
where $N\ge 2$. These groups are modular and their index is given by
\[ |\Gamma(1):\Gamma(N)|=  \{ \begin{array}{cc}6& \mbox{if }N=2\\ 1/2N^3
\prod_{p|N}(1-1/p^{2})& \mbox{if }N>2 \end{array} \]
\[ |\Gamma(1):\Gamma_{0}(N)|=
N \prod_{p|N}(1+1/p) \]
where the product runs over all prime divisors $p$ of $N$.\\
Later on we will need the generators for the groups
$\Gamma(1),~\Gamma(2)$ and $\Gamma_{0}(2)$.  The group $\Gamma(1)$ is
generated by the elements
\[ Q= \left( \begin{array}{cc}0&1 \\ -1 & 0\end{array}  \right)\mbox{ and }
T=\left( \begin{array}{cc}1 & 1 \\ 0 & 1 \end{array} \right) \]
fulfilling the relations \[ Q^2=(QT)^3=I. \] The group $\Gamma(2)$ is
freely generated by
\[ A_{1}=T^{2}\quad \mbox{ and } \quad A_{2}=QT^{2}Q \] 
and $\Gamma_{0}(2)$ is generated by
\[ B_{1}=T\quad,\quad B_{2}=QT^{-2}Q\quad \mbox{and}\quad B_{3}=T^{-1}B_{2}\]
which fulfill the relations
\[ B_{1}B_{2}B_{3}=I \quad \mbox{ and } \quad B_{3}^{2}=I .\]
Now let us mention a few facts about the coset sets.  Consider the map
\[ Mod_{N}:PSL(2,\Z)\longmapsto PSL(2,\Z/N\Z) \mbox{ with }
Mod_{N}(A)=A \mbox { mod } N\] This is a surjective homomorphism with
kernel $\Gamma(N)$. Hence $\Gamma(N)$ is a normal subgroup of
$\Gamma(1)$
and the group $\tilde\Gamma(N)=\Gamma(N)\backslash\Gamma(1)$ is isomorphic to $PSL(2,\Z/N\Z)$.\\
Now consider $\Gamma_{0}(N)$. Obviously we have
$\Gamma(N)\subseteq\Gamma_{0}(N)\subseteq\Gamma(1)$ and
$Mod_{N}(\Gamma_{0}(N))=\Theta(N)$ where $\Theta(N)$ is the group
\[  \{ \left( \begin{array}{cc}a&b \\ 0& a^{-1}\end{array}  \right)|~b\in\Z/N\Z,\quad a \in(Z/N\Z)^{\star}\}. \]
$\Gamma_{0}(N)$ is not normal in $\Gamma(1)$ but there is a natural
bijection from $\tilde\Gamma_{0}(N)=\Gamma_{0}(N)\backslash\Gamma(1)$
to $\Theta(N)\backslash\tilde\Gamma(1)$. For $N=p$ prime a system of
representatives of this coset set is given by
\[ \{QT^{i}|i=0\dots p-1\}\cup\{I\}. \]
\\
To each modular group $\Gamma$ there corresponds a modular surface
given by the quotient space $\Gamma\backslash \H$ which consists of
the equivalence classes of $\Gamma$-equivalent points of $\H$. This
surface is topologically a sphere with finitely many handles and
finitely many cusps. The cusps are located at the rationals on the
real axis and at $\infty$.  The compactification of this surface is a
Riemannian surface with the Riemannian metric coming from the Poincare
metric on $\H$.  Let $g$ be the genus of this surface, $v_{\infty}$
the number of cusps, $v_{2}$ be the number of elliptic elements of order
2 in $\Gamma$ and $v_{3}$ be the number of elliptic elements of order $3$
in $\Gamma$\footnote{An element in $\Gamma$ is called elliptic if it
  has the fixed points $z$ and $\bar z$ with $z\in\H$.}. With these
notations we have \cite{[MI]}
\[ g=1+\frac{|\Gamma(1):\Gamma|}{12}-\frac{v_{2}}{4}-\frac{v_{3}}{3}-\frac{v_{\infty}}{2}. \]
In the case of the principal congruence subgroups and congruence
subgroups of Hecke type all these quantities are explicitly known. For
$\Gamma(N)$ we have $v_{2}=v_{3}=0$ and
$v_{\infty}=|\Gamma(1):\Gamma(N))|/N$ and hence
\[ g(\Gamma(N))=1+\frac{|\Gamma(1):\Gamma(N)|}{12}-\frac{|\Gamma(1):\Gamma(N)| }{2N}. \] 
For $\Gamma_{0}(N)$ we have
\[ v_{2}(\Gamma_{0}(N))=  \{ \begin{array}{cc}0& \mbox{if }4|N\\ 
\prod_{p|N}(1+\{\frac{-1}{p}\})& \mbox{if }4 \not\vert N\end{array} \]
\[ v_{3}(\Gamma_{0}(N))=  \{ \begin{array}{cc}0& \mbox{if }9|N\\ 
\prod_{p|N}(1+\{\frac{-3}{p}\})& \mbox{if }9\not\vert N \end{array} \]
and
\[ v_{\infty}(\Gamma_{0}(N))=\sum_{0<d|N}\phi(d,N/d)\]
where $\{-\}$ denotes here the quadratic residue symbol and $\phi$
denotes the Euler function. Again the genus of the surfaces
corresponding to $\Gamma_{0}(N)$ can be explicitly calculated using
these expressions.
\section{Automorphic forms and the Eichler cohomology for modular groups}
Let $\Gamma$ be a modular group.  We will first recall the definition
of {\bf automorphic forms} and {\bf cusp forms} on $\Gamma$. For
$A=\left(\begin{array}{cc} a&b \\c&d\end{array}\right)\in\Gamma$,
$f:\H\longmapsto \C$ and $k\in\Z$ we set
\[ [A]_{k}f(z)=(cz+d)^{-k}f(Az).\]
$f$ is called an automorphic form of weight $k$ (or degree $-k$) for
$\Gamma$ if $f$ is holomorphic on $\H$ and at the cusps of
$\Gamma\backslash\H$ such that $[A]_{k}f=f$ for all $A\in \Gamma$. If
in addition $f$ vanishes at the cusps of $\Gamma$ it is called cusp
form.  We denote by $A_{k}(\Gamma)$ the space of automorphic forms and
by $C_{k}(\Gamma)$ the space of cusp forms.  It is possible to
calculate the dimension of these spaces by relating them to
differentials on $\Gamma\backslash\H$ and using the Riemann-Roch
Theorem. In fact we have the following Proposition, see \cite{[MI]}.
\begin{proposition}
  For arbitrary modular groups $\Gamma$ we have with the above
  notations
\[ \dim C_{k}(\Gamma)=(k-1)(g-1)+(k/2-1)v_{\infty}+[k/4]v_{2}+[k/3]v_{3}\]
for all $k> 2$ even and $\dim C_{2}(\Gamma)=g$. Moreover
\[  \dim A_{k}(\Gamma)=\dim C_{k}(\Gamma)+v_{\infty}(\Gamma) \]
and the automorphic forms that are not cusps forms are given by
Eisenstein series.
\end{proposition}
For the principal congruence subgroups $\Gamma(N)$ we get the
following corollary.
\begin{corollary}
  For all $N\ge 2$ and all $k\in \N$ even we have
\[ \dim C_{k}(\Gamma(N))=\frac{|\Gamma(1):\Gamma(N)|}{12}(k-1)-\frac{v_{\infty}(\Gamma(N))}{2}. \]
\end{corollary}
Using Proposition 3.1 it is also possible to calculate $\dim
C_{k}(\Gamma_{0}(N))$.  In the appendix of \cite{[MI]} the reader will
find a table with the explicit values of
$\dim C_{k}(\Gamma_{0}(N))$ for a wide range of $N$ and $k$.\\

Now we introduce the {\bf Eichler cohomology}, see \cite{[EI]}.  Let
$f\in C_{k+2}(\Gamma)$ be a cusp form of weight $k+2$ (degree
$-(k+2)$) with $k\in\N_{0}$ an even number.  We are interested in the
$k+1$-th undetermined integral of $f$. Obviously the integral
\[\Theta(\tau)=\frac{1}{k!}\int_{\tau_{0}}^{\tau}(t-z)^{k}f(z)dz\]
is path independent for $\tau,\tau_{0}\in\H$ and
\[\frac{d^{k+1}}{d\tau^{k+1}}\Theta(\tau)=f(\tau). \]
Let $\P_{k}$ denote the linear space of all polynomials over $\C$ of
degree less or equal $k$. We can add an arbitrary polynomial
$\phi\in\P_{k}$ to the above integral as an ``integration constant''
to obtain an arbitrary $k+1$-th integral of $f$
\[\Theta(\tau)=\frac{1}{k!}\int_{\tau_{0}}^{\tau}(t-z)^{k}f(z)dz+\phi(\tau).\]
Moreover it is easy to show that
\[ \Omega_{A}(\tau):=[A]_{-k}\Theta(\tau)-\Theta(\tau) \]
is a polynomial in $\P_{k}$ for all $A\in\Gamma$.  The map
\[ \Omega:A\longmapsto \Omega_{A}\]
from $\Gamma$ into $\P_{k}$ is a cocycle since
\[  \Omega_{AB}=[B]_{-k}\Omega_{A}+\Omega_{B}. \]
We denote the space of all cocycles by $\bar E_{k}(\Gamma)$;
\[\bar E_{k}(\Gamma)=\{\Omega:\Gamma\longmapsto \P_{k}|\Omega_{AB}=[B]_{-k}\Omega_{A}+\Omega_{B}\}\]
Special cocycles are the coboundaries
\[ \Omega_{A}(\tau)=[A]_{-k}\phi(\tau)-\phi(\tau) \]
for $\phi\in\P_{k}$.  The space of all Eichler cohomology classes is
now the quotient
\[ E_{k}(\Gamma):=\bar E_{k}(\Gamma)/\P_{k}\]
where we identify $\P_{k}$ with the space of coboundaries in the
obvious way. Now every cusp form $f\in C_{k+2}(\Gamma)$ corresponds to
a cohomology class in $E_{k}(\Gamma)$.  Eichler's work \cite{[EI]}
shows that this correspondence is two to one.  Moreover \cite{[EI]}
implicitly contains the following result.
\begin{proposition}
  For all $k\ge 2$ even we have
\[ \dim E_{k}(\Gamma)=2 \dim C_{k+2}(\Gamma)+ v_{\infty}(\Gamma) \]
and especially
\[ \dim E_{k}(\Gamma(N))=\frac{|\Gamma(1):\Gamma(N)|}{6}(k+1)\]
\end{proposition}
For $\Gamma(2)$ and $\Gamma_{0}(2)$ we want to compute the spaces
$E_{k}$ explicitly.  Let $A_{1}$ and $A_{2}$ be the generators of
$\Gamma(2)$. Since $\Gamma(2)$ is freely generated by these matrices
it is obvious that every pair of polynomials
$\Omega_{A_{1}},\Omega_{A_{2}}\in\P_{k}$ determines a cocycle and vice
versa. Hence
\[ E_{k}(\Gamma(2))=\bar E_{k}(\Gamma(2))/\P_{k}\cong (\P_{k}\times\P_{k})/\P_{k}\cong \P_{k}. \] 
Let $B_{1}, B_{2}, B_{3}$ be the generators of $\Gamma_{0}(2)$.
Choose arbitrary polynomials $\Omega_{B_{1}}, \Omega_{B_{3}} \in
P_{k}$.  By the relation $B_{1}B_{2}B_{3}=I$ a polynomial
$\Omega_{B_{2}}\in \P_{k}$ is uniquely determined by the cocycle
condition. By the relation $B_{3}^{2}=I$, $\Omega_{B_{3}}$ fulfills
the cocycle relation if and only if
$[B_{3}]_{-k}\Omega_{B_{3}}+\Omega_{B_{3}}=0$ which means
\[ \Omega_{B_{3}}(\frac{z+1}{-2z-1})(-2z-1)^{k}+\Omega_{B_{3}}(z)=0.\]
Let $\Upsilon_{k}$ be the linear space of all polynomials
$\Omega\in\P_{k}$ obeying this functional equation. We now see that a
cocycle for $\Gamma_{0}(2)$ is uniquely determined by an arbitrary
polynomial in $P_{k}$ and an polynomial in $\Upsilon_{k}$. Hence
\[ E_{k}(\Gamma_{0}(2))=\bar E_{k}(\Gamma_{0}(2))/\P_{k}\cong (\P_{k}\times\Upsilon_{k})/\P_{k}\cong \Upsilon_{k} .\] 
By Proposition 3.1 and 3.2 we now get a formula for the dimension of
the space $\Upsilon_{k}$
\[ \dim\Upsilon_{k}=2[(k-1)/4]+2. \]
\section{The Transfer operator for modular groups}
Here we consider the transfer operator for modular groups $\Gamma$
introduced by Chang and Mayer in \cite{[CM1]}, \cite{[CM2]} and
\cite{[CM3]}. In these papers the
reader will find a detailed discussion  of the properties of the transfer operator that are mentioned here.\\
For a function $f:\C\times\tilde\Gamma\times\Z/2\Z\longmapsto \C$ and
$s\in\C$ we define a formal operator $L_{s}$ by
\[ L_{s}f(z,A,\epsilon)=\sum_{n=1}^{\infty}(\frac{1}{z+n})^{2s}f(\frac{1}{z+n},AQT^{n\epsilon},-\epsilon)\]
where $\tilde\Gamma=\Gamma\backslash\Gamma(1)=\{A_{1},A_{2},\dots,A_{\mu}\}$ with $\mu=|\Gamma(1):\Gamma|$.
We call this operator the {\bf transfer operator}. There is another
way to express this operator using representation theory\footnote{We
  refer to \cite{[FH]} for an introduction to representation theory.}.
We define $\chi^{\Gamma}$ as the representation of $\Gamma(1)$ which is
induced by the trivial representation of $\Gamma$. That means
$ \chi^{\Gamma}:\Gamma(1)\longmapsto \{0,1\}^{\mu}$ is given by
\[ \chi^{\Gamma}(G)=(\chi(A_{i}G A_{j}))_{i,j=1,\dots,\mu} \]
with 
\[ \chi(G)=\{ \begin{array}{cc} 1\quad G\in \Gamma\\ 0\quad G\not\in\Gamma \end{array}.  \]
We can now write
the operator $L_{s}$ in the form
\[  L_{s}f(z,\epsilon)=\sum_{n=1}^{\infty}(\frac{1}{z+n})^{2s}\chi^{\Gamma}(QT^{n\epsilon})f(\frac{1}{z+n},-\epsilon)\]
which acts on functions $f:\C\times\Z/2\Z\longmapsto\C^{\mu}$.  Now let
\[ D=\{z|z\in\C,|z-1|<3/2 \}\]
and consider the Banach space
\[\tilde B:= B(D\times \Z/2\Z)^{\mu}\]
where $B(D\times \Z/2\Z)$ is the space of all complex valued functions
holomorphic on $D\times \Z/2\Z$ and continuous on the boundary of this
set.  The operator $L_{s}$ is well defined on $\tilde B$ and
holomorphic if $\Re(\beta)>1/2$. It can be continued to an operator
which is meromorphic on the whole complex $\beta$-plane with possible
poles only for $\beta=\beta(\kappa)=(1-\kappa)/2$ with
$\kappa\in\N_{0}$.  Moreover it was shown by Chang and Mayer that
$L_{s}$ is a nuclear operator and hence of trace class. The Fredholm
determinant of the operator is given
by Selberg's zeta function with respect to the representation $\chi^{\Gamma}$, see \cite{[VE]}.  \\
The transfer operator has a well known interpretation within the
theory of dynamical systems. Let us consider the geodesic flow on the
surface $\Gamma\backslash \H$.  This is a Hamiltonian flow and
provides a classical example of a ``chaotic'' dynamical systems, see
\cite{[KH]}.  It is possible to take a Poincare section for the flow
in order to obtain a discrete time system defined by the first return
map.  In appropriate coordinates this map is given by
\[ g(z,A,\epsilon)=(\frac{1}{z}-[\frac{1}{z}],AQT^{n\epsilon},-\epsilon ) \]
on the space $I_{2}\times\tilde\Gamma\times\Z/2\Z$. Now the Transfer
operator defined above is nothing but a generalization of the
classical Perron Frobenius operator in ergodic theory for this system, see
\cite{[KH]} and \cite{[LM]}.
\section{Eigenfunctions of the transfer operator and the Lewis equation}
The key tool for determining eigenfunctions of the transfer operator
for a fixed parameter $s \in \C$ to an eigenvalue $\lambda=\lambda(s)$
is a certain functional equation, which is called {\bf Lewis
  equation}. For $\Gamma(1)$ this equation was found in \cite{[LW]}.
For modular groups $\Gamma$ we have the following Proposition
contained in \cite{[CM3]} which relates eigenfunctions of the transfer
operator to solutions of this Lewis equation.
\begin{proposition}
  $f\in\tilde B$ is a eigenfunction of the transfer operator $L_{s}$
  with eigenvalue $\lambda$ i.e. $L_{s}f=\lambda f$ if and only if
\[\lambda( f( z,\epsilon)- \chi^{\Gamma}(QT^{\epsilon}Q)f(z+1,\epsilon))=(\frac{1}{z+1})^{2s}\chi^{\Gamma}(QT^{\epsilon})f(\frac{1}{z+1},-\epsilon) \]
for $\epsilon\in \Z/2\Z$ and
\[\lim_{z\longmapsto \infty}(\lambda f(z,\epsilon)-\sum_{l=0}^{\kappa}\sum_{m=1}^{r}(1/r)^{2s+l}\chi^{\Gamma}(QT^{m\epsilon})\frac{f^{(l)}(0,-\epsilon)}{l!}\zeta(2s+l,\frac{z+m}{r}))=0\]
for $\epsilon\in \Z/2\Z$ and $\kappa >2\Re(s)$ where $\zeta$ is the
Hurwitz zeta function.
\end{proposition}
{\bf Remark 5.1~}
  In the special case $\chi^{\Gamma}(T^{2})=I$ the two equations in
  the last Proposition are the same for $\epsilon=1$ and $\epsilon=-1$
  and thus reduce to
\[\lambda( f( z)- \chi^{\Gamma}(QTQ)f(z+1))=(\frac{1}{z+1})^{2s}\chi^{\Gamma}(QT)f(\frac{1}{z+1}) \]
and
\[\lim_{z\longmapsto \infty}(\lambda f(z)-\sum_{l=0}^{\kappa}\sum_{m=1}^{r}(1/r)^{2s+l}\chi^{\Gamma}(QT^{m})\frac{f^{(l)}(0)}{l!}\zeta(2s+l,\frac{z+m}{r}))=0
~\kappa>2\Re(s).\] For the groups $\Gamma(2)$ and $\Gamma_{0}(2)$ we
indeed have $\chi^{\Gamma}(T^{2})=I$ since $T^{2}=I$ in
$\tilde\Gamma(2)$ and $\tilde\Gamma_{0}(2)$. We will use this fact
later on.\\\\
{\bf Remark 5.2}
  Consider the following functional equation
\[ g(z)-\chi^{\Gamma}(QTQ)g(z+1)=(\frac{1}{z+2})^{2s}\chi^{\Gamma}(QT^{2})g(\frac{-1}{z+2}). \]
Following \cite{[LZ]} we call this equation the {\bf Master equation}.
Let $g$ be a solution of the Master equation and set
\[g^{+}(z,\epsilon)=\{ \begin{array}{cc} g(z)& \mbox{if }\epsilon=1 \\ (1/(z+1))^{2s}\chi^{\Gamma}(T)g(-z/(z+1))&\mbox{ if }\epsilon=-1  \end{array} \]
and
\[ g^{-}(z,\epsilon)=\{ \begin{array}{cc} g(z)& \mbox{if }\epsilon=-1 \\ -(1/(z+1))^{2s}\chi^{\Gamma}(T)g(-z/(z+1))&\mbox{ if }\epsilon=-1  \end{array} \]
Then $g^{+}$ is a solution of the Lewis equation for $\lambda=+1$ and
$g^{-}$ is a solution of the Lewis equation for $\lambda=-1$. In fact
also the converse is true. If $g^{+}$ is a solution of the Lewis
equation with $\lambda=1$ then $g(z):=g^{+}(z,1)$ is a solution of the
Master equation. The same is true for $g^{-}$. See section 3 of
\cite{[CM3]}.
\section{Polynomial Eigenfunctions of the transfer operator}
We are interested in this paper primarily in polynomial eigenfunctions
of the transfer operator. By Proposition 5.1 we see that the transfer
operator $L_{s}$ may have eigenfunctions in $\P_{n}\times\{-1,1\}$ if
and only if $s=-n/2$ with $n\in \N$ even. Moreover for all polynomial
eigenfunctions
the asymptotic condition in Proposition 5.1~is trivially satisfied.\\
By well known properties of Selberg's zeta function one knows that
$L_{-n/2}$ has the eigenvalue $+1$ and $-1$, see \cite{[CM1]}. We
restrict our attention from now on to these eigenvalues.  By
$\wp^{+}_{n}(\Gamma)$ resp. $\wp^{-}_{n}(\Gamma)$ we denote the space
of all polynomial eigenfunctions of degree $n$ of the transfer
operator $L_{-n/2}$ for the modular group $\Gamma$ to the eigenvalue
$+1$ resp. $-1$. More precisely
\[ \wp_{n}^{+}(\Gamma)=\{p\in (\P_{n})^\mu|L_{-n/2}p^{+}=p^{+}\}\]
\[ \wp_{n}^{-}(\Gamma)=\{p\in (\P_{n})^\mu|L_{-n/2}p^{-}=-p^{-}\}\]       
where $\mu=|\Gamma(1):\Gamma|$. Furthermore let
\[ \wp_{n}(\Gamma)= \wp_{n}^{+}(\Gamma)\bigoplus\wp_{n}^{-}(\Gamma).\]
In the next sections we will study this space for principal congruence
subgroups and congruence subgroups of Hecke type.
\section{The case of principal congruence subgroups}
In this section we consider the principal congruence group $\Gamma(N)$
defined in section 2 for $N\ge 2$.  We are interested in the dimension
of the space $\wp_{n}(\Gamma(N))$ of all polynomial eigenfunctions of
degree $n$ of the transfer operator $L_{-n/2}$ where $n$ is an even
number (see section 6) and in the relation of this space to the space
of Eichler's cohomology classes $E_{n}(\Gamma(N))$ (see section 4).
We have the following conjecture.
\begin{conjecture} For all $N\ge2$ and $n\ge 2$ even we have
\[ \dim \wp_{n}(\Gamma(N))=\frac{|\Gamma(1):\Gamma(N)|}{6}(n+1)\]
and hence
\[ \wp_{n}(\Gamma(N))\cong E_{n}(\Gamma(N). \]
\end{conjecture}
{\bf The case $N=2$}\\\\
In the case $N=2$ we can prove this conjecture. Since
$|\Gamma(1):\Gamma(2)|=6$ it is in fact an immediate consequence of
the following Proposition
\begin{proposition}
  We have
\[  \dim \wp_{n}^{-}(\Gamma(2))=n/2+2 \quad\mbox{ and }\quad \dim \wp_{n}^{+}(\Gamma(2))=n/2-1. \]
\end{proposition}
\proof  Using the fact that $\tilde\Gamma(2)$ is isomorphic to
$PSL(2,\Z/2\Z)$ it is easy to calculate the induced representation
$\chi^{\Gamma(2)}$. We get
\[ \chi^{\Gamma(2)}(QTQ)=\left( \begin{array}{cccccc}0&0&0&0&0&1\\0&0&0&0&1&0\\0&0&0&1&0&0\\0&0&1&0&0&0 \\ 0&1&0&0&0&0 \\ 1&0&0&0&0&0  \end{array} \right) \mbox{ and }\chi^{\Gamma(2)}(QT)=\left( \begin{array}{cccccc}0&0&0&0&1&0\\0&0&0&0&0&1\\0&1&0&0&0&0\\1&0&0&0&0&0 \\ 0&0&0&1&0&0 \\ 0&0&1&0&0&0  \end{array} \right)\]
If we set $f(z)=\phi(z+1)$ the Lewis equation from remark 5.1 is given
by the following system of functional equations
\[\begin{array}{cccccc}\lambda(\phi_{1}(z)-\phi_{6}(z+1))-z^{n}\phi_{5}(\frac{z+1}{z})=0\\ \lambda(\phi_{2}(z)-\phi_{5}(z+1))-z^{n}\phi_{6}(\frac{z+1}{z})=0\\ \lambda(\phi_{3}(z)-\phi_{4}(z+1))-z^{n}\phi_{2}(\frac{z+1}{z})=0\\ \lambda(\phi_{4}(z)-\phi_{3}(z+1))-z^{n}\phi_{1}(\frac{z+1}{z})=0\\ \lambda(\phi_{5}(z)-\phi_{2}(z+1))-z^{n}\phi_{4}(\frac{z+1}{z})=0\\ \lambda(\phi_{6}(z)-\phi_{1}(z+1))-z^{n}\phi_{3}(\frac{z+1}{z})=0.  \end{array} \]
First note that that $\phi_{1}$ and $\phi_{2}$ are uniquely determined
by the first two equations if $\phi_{5}$ and $\phi_{6}$ are given.  By
substituting $1/z$ for $z$ and multiplying with $-\lambda z^{n}$ the
last two equation are equivalent to
\[\begin{array}{cc}z^{n}(-\phi_{5}(\frac{1}{z})+\phi_{2}(\frac{z+1}{z}))+\lambda\phi_{4}(z+1)=0\\ z^{n}(-\phi_{6}(\frac{1}{z})+\phi_{1}(\frac{z+1}{z}))+\lambda\phi_{3}(z+1)=0. \end{array} \]
Now adding the third and the fourth equation to this equation we see
that
\[\begin{array}{cc} \phi_{3}(z)=\lambda z^{n}\phi_{5}(\frac{1}{z}) \\ \phi_{4}(z)=\lambda z^{n}\phi_{6}(\frac{1}{z}).\\  \end{array} \]
Thus $\phi_{3}$ and $\phi_{4}$ are uniquely determined if $\phi_{5}$
and $\phi_{6}$ are given.  Now inserting the expressions for
$\phi_{1},\phi_{2},\phi_{3}$ and $\phi_{4}$ in the last two equations
we get
\[\begin{array}{cc}\phi_{5}(z)-\phi_{5}(z+2)=(z+1)^{n}(\lambda\phi_{6}(\frac{z+2}{z+1})+\phi_{6}(\frac{z}{z+1}))\\ \phi_{6}(z)-\phi_{6}(z+2)=(z+1)^{n}(\lambda\phi_{5}(\frac{z+2}{z+1})+\phi_{5}(\frac{z}{z+1}))  \end{array} \] 
and substituting $z-1$ for $z$ this gives
\[\begin{array}{cc}\phi_{5}(z-1)-\phi_{5}(z+1)=z^{n}(\lambda\phi_{6}(\frac{z+1}{z})+\phi_{6}(\frac{z-1}{z}))\\ \phi_{6}(z-1)-\phi_{6}(z+1)=z^{n}(\lambda\phi_{5}(\frac{z+1}{z})+\phi_{5}(\frac{z-1}{z}))  \end{array} \]
With
\[\bar T\phi(z):=\phi(z-1)-\phi(z+1)  \quad\mbox{ and }\quad G\phi(z):=z^{n}(\lambda\phi(\frac{z+1}{z})+\phi(\frac{z-1}{z})).\]we can write this system in the form
\[ \bar T\phi_{5}= G\phi_{6}\quad\mbox{ and }\quad \bar T\phi_{6}= G\phi_{5} \] 
We have to determine the solutions
$(\phi_{5},\phi_{6})\in \P_{n}^{2}$ of this system of functional equations.\\
Obviously $G$ is a linear map of $\P_{n}$ into $\P_{n}$ and $\bar
T$ is a linear map $\P_{n}$ onto $\P_{n-1}$. Let
$\P_{n}^{e}=\{\phi\in\P_{n}|\phi\mbox{ is even}\}$ and
$\P_{n}^{o}=\{\phi\in\P_{n}|\phi\mbox{ is odd}\}$. If
$\phi\in\P^{e}_{n}$ we have
\[ \bar T\phi(-z)=\phi(-z-1)-\phi(-z+1)=\phi(z+1)-\phi(z-1)=-\bar T\phi(z)\]
and hence $\bar T\phi\in\P^{o}_{n-1}$ and if $\phi\in\P^{o}_{n}$ we
have
\[ \bar T\phi(-z)=\phi(-z-1)-\phi(-z+1)=-\phi(z+1)+\phi(z-1)=\bar T\phi(z)\]
and hence $\bar T\phi\in\P^{e}_{n-1}$.\\
Now consider the case $\lambda=-1$. In this case $G\phi\in\P_{n-1}$
and since
\[ G\phi(-z)=(-z)^{n}(-\phi(\frac{-z+1}{-z})+\phi(\frac{-z-1}{-z}))\]\[=z^{n}(\phi(\frac{z+1}{z})-\phi(\frac{z-1}{z}))=-G\phi(z)\]
the map $G$ is onto $\P^{o}_{n-1}$. \\
Let $(\phi_{5},\phi_{6})\in \P_{n}^{2}$ be solutions of $\bar
T\phi_{5}= G\phi_{6}$ and $\bar T\phi_{6}= G\phi_{5}$. We know that
$G\phi_{6}$ and hence $\bar T\phi_{5}$ is odd. Now assume that
$\phi_{5}=\phi_{5}^{e}+\phi_{5}^{o}$ where $\phi_{5}^{e}$ is an even
polynomial and $\phi_{5}^{o}$ is an odd polynomial and not zero. We
have $\bar T\phi_{5}=\bar T\phi_{5}^{e}+\bar T\phi_{5}^{o}$ where
$\bar T\phi_{5}^{e}$ is odd but $\bar T\phi_{5}^{o}$ is even and not
zero.  Hence $\bar T\phi_{5}$ would not be odd. This is a
contradiction and hence $\phi_{5}$ has to be even.
By the same argument we can show that $\phi_{6}^{e}$ is even.\\
Now let $\phi_{5}$ be an arbitrary polynomial in $\P^{e}_{n}$.  We
have $\bar T\phi_{5}\in \P^{o}_{n-1}$ and the map
$G:\P_{n}^{e}\longmapsto P_{n-1}^{o}$ is onto and has a kernel
consisting of all constant polynomials. Hence we see that all
solutions $\phi_{6}\in\P_{n}$ of $\bar T\phi_{5}= G\phi_{6}$ are given
by $\phi_{6}=\bar\phi_{6}+d$ where $\bar\phi_{6}\in \P_{n}^{e}$ is
uniquely determined and $d$ is an arbitrary constant. Furthermore if
$\bar T\phi_{5}= G\phi_{6}$ we have
\[ \phi_{5}(z-1)-\phi_{5}(z+1)=z^{n}(-\phi_{6}(\frac{z+1}{z})+\phi_{6}(\frac{z-1}{z})).\]
Substituting $1/z$ for $z$ and using the fact that $\phi_{5}$ and
$\phi_{6}$ are even this implies
\[ \phi_{6}(z-1)-\phi_{6}(z+1)=z^{n}(\phi_{5}(-\frac{z+1}{z})+\phi_{5}(\frac{z-1}{z})). \]
This means that the functional equation $\bar T\phi_{6}= G\phi_{5}$
holds as well. Hence the space of solutions in $\P_{n}$ of our system
of functional equations is isomorphic to $\P_{n}^{e}\times \P_{o}$ and
hence $\wp_{n}(\Gamma(2))\cong\P_{n}^{e}\times \P_{o}$.
Since $\dim\P_{n}^{e}=n/2+1$ for $n$ even this completes the proof of our  result in the case $\lambda=-1$. \\
In the case $\lambda=+1$ an argument along the same lines
just exchanging the role of odd and even polynomials gives the desired result. \qed\\\\
{\bf The case $N>2$}\\\\
In the case $N>2$ we have some numerical evidence for our general
conjecture.  Table 1 below contains the dimension of
$\wp_{n}(\Gamma(N))$ for some values of $N>2$ and $n$ even.  By the
formula for $|\Gamma(1):\Gamma(N)|$ given in section 2 we see from
table 1 that conjecture 7.1 is really true for all values
of $N$ and $n$ we have checked.\\
\begin{center}
\begin{tabular}{|c|c|c|c|c|c|c|}
\hline
      ~&  $n=2$ & $n=4$ & $n=6$ & $n=8$ & $n=10$ & $n=12$    \\      
\hline
N=3 & 6  & 10 & 14 & 18 & 22          & 24     \\ 
\hline
N=4 & 12 & 20 & 28 & 36 & 44           & 52    \\ 
\hline
N=5 &  30 & 50  & 70 & 90 &  110           & 130   \\ 
\hline
N=6 &  36 & 60  & 84 & 108 &  132      & 156  \\ 
\hline
N=7 &  84 & 140  &196  &  242   &               &  \\ 
\hline
N=8 &  96 & 160  &  &    &               &  \\ 
\hline
\end{tabular}
\end{center}
{\bf Table 1: $\dim\wp_{n}(\Gamma(N))$}\\ \\
We have calculated these dimensions algebraically using Matematica.
Let us make a few comments on this computation.  It is easy to compute
the finite group $PSL(2,\Z/N\Z)$ and the action of $Q$ and $T$ on
these group. This gives us the regular representation of
$PSL(2,\Z/N\Z)$ which can be identified with the induced
representation $\chi^{\Gamma}$. Now we can explicitly determine the
Lewis equation for $\Gamma(N)$. Inserting polynomials into this
equation we get a system of $C(N,n)=|\Gamma(1):\Gamma(N)|n\approx
N^{3}n$ linear equation. This system can be solved algebraically in at
most a day on a PC if
$C(N,n)$ is less then  $2000$.\\
\section{The case of congruence subgroups of Hecke type}
In this section we consider the congruence groups of Hecke type
$\Gamma_{0}(N)$ defined in section 2 for $N\ge 2$.  We are interested
in the relation of the space $\wp_{n}(\Gamma_{0}(N))$ of all
polynomial eigenfunctions of degree $n$ of the transfer operator
$L_{-n/2}$ where $n$ is an even number to the space of Eichler's
cohomology classes $E_{n}(\Gamma_{0}(N))$, see sections 4 and 6.  We
have the following conjecture.
\begin{conjecture} For all $N\ge2$ and $n\ge 2$ even we have
\[ \wp_{n}(\Gamma_{0}(N))\cong E_{n}(\Gamma_{0}(N)). \]
\end{conjecture}
{\bf The case $N=2$}\\\\
Let
\[ \tilde\wp_{n}^{+}=\{p\in (\P_{n+1})^\mu|L_{-n/2}p^{+}=p^{+}\}\]
\[ \tilde\wp_{n}^{-}=\{p\in (\P_{n+1})^\mu|L_{-n/2}p^{-}=-p^{-}\}\]   
where $L$ is the transfer operator with respect to $\Gamma_{0}(2)$ and
$p^{+}$ and $p^{-}$ are defined in remark 5.2~.  We consider here
polynomials in $\P_{n+1}$ instead of polynomials in $\P_{n}$ for some
technical reasons.  In fact there exist such eigenfunctions of
$L_{-n/2}$ and we can explicitly relate the space $ \tilde\wp_{n}(N)=
\tilde\wp_{n}^{+}(N)\oplus \tilde\wp_{n}^{-}(N)$ to the space
$\Upsilon_{n}$ which is isomorphic to Eichler's cohomology classes for
$\Gamma_{0}(2)$, see section 3.
\begin{proposition}
  Let $n\in\N$ be an even number and let $\tilde T$ be the linear
  operator given by $\tilde T\phi(z)=\phi(z)-\phi(z+1)$ from
  $\P_{n+1}$ to $\P_{n}$. Then
\[ \Upsilon_{n}\cong \tilde T \tilde\wp_{n}. \]
\end{proposition}
\proof Using the fact that a system of representatives of
$\tilde\Gamma_{0}(2)$ is given by $\{I,Q,QT \}$ it is easy to
calculate the induced representation $\chi^{\Gamma(2)}$. We get
\[ \chi^{\Gamma(2)}(QTQ)=\left( \begin{array}{ccc}0&0&1\\0&1&0\\1&0&0\\ \end{array} \right) \mbox{ and }\chi^{\Gamma(2)}(QT)=\left( \begin{array}{ccc}0&0&1\\1&0&0\\0&1&0 \end{array} \right). \]
If we set $f(z)=\phi(z+1)$ the Lewis equation from remark 5.1 is given
by the following system of functional equations
\[\begin{array}{ccc}\lambda(\phi_{1}(z)-\phi_{3}(z+1))-z^{n}\phi_{3}(\frac{z+1}{z})=0\\ \lambda(\phi_{2}(z)-\phi_{2}(z+1))-z^{n}\phi_{1}(\frac{z+1}{z})=0\\ \lambda(\phi_{3}(z)-\phi_{1}(z+1))-z^{n}\phi_{2}(\frac{z+1}{z})=0.\\  \end{array} \]
Expressing $\phi_{3}$ through $\phi_{2}$ and $\phi_{1}$ using the
third equation and inserting $\phi_{3}$ into the second equation
yields
\[ \phi_{3}(z)=\lambda z^{n}\phi_{2}(1/z). \]
Thus $\phi_{3}$ can be expressed through $\phi_{2}$. Moreover by the
second equation $\phi_{1}$ can be expressed through $\phi_{2}$,
\[ \phi_{1}(z)=\lambda(z-1)^{n}(\phi_{2}(1/(z-1))-\phi_{2}(z/(z-1)))\]
Inserting this expression into the first equation yields a functional
equation for $\phi_{2}$
\[ \phi_{2}(z)-\phi_{2}(z+1)=(2z+1)^{n}(\lambda\phi_{2}(\frac{z+1}{2z+1})+\phi_{2}(\frac{z}{2z+1}))\]
This shows
\[ \tilde\wp_{n}^{+}\cong\{\phi\in\P_{n+1}|\phi(z)-\phi(z+1)=(2z+1)^{n}(\phi(\frac{z+1}{2z+1})+\phi(\frac{z}{2z+1}))\}\]
and
\[  \tilde\wp_{n}^{-}\cong\{\phi\in\P_{n+1}|\phi(z)-\phi(z+1)=(2z+1)^{n}(-\phi(\frac{z+1}{2z+1})+\phi(\frac{z}{2z+1}))\}                   \]
We will identify these isomorphic spaces in the following. We want to
show that
\[ \tilde\wp_{n}^{+}\subseteq \P_{n+1}^{o} \mbox{ and } \tilde\wp_{n}^{-}\subseteq \P_{n+1}^{e} \] 
Let $\phi\in\tilde\wp^{+}_{n}$. Substituting $-z-1$ in the functional
equation determining $\phi$ and using the fact that $n$ is even we get
\[\phi(-z-1)-\phi(-z)=(2z+1)^{n}(\phi(\frac{z+1}{2z+1})+\phi(\frac{z}{2z+1})) \]
and hence
\[ \phi(-z-1)-\phi(-z)=\phi(z)-\phi(z+1)=-(\phi(z+1)-\phi(z)).\]
Decomposing $\phi$ in an even and an odd polynomial we see that this
implies $\phi=\phi^{o}+c$ where $\phi^{o}$ is odd and $c$ is a
constant.
Inserting this into the functional equitation we get $2(2z+1)^{n}c=0$. Hence $c=0$ and $\phi\in\P^{o}_{n+1}$.\\
Let $\phi\in\tilde\wp^{-}_{n}$. Substituting $-z-1$ in the functional
equation determining $\phi$ and using the fact that $n$ is even we get
\[\phi(-z-1)-\phi(-z)=(2z+1)^{n}(\phi(\frac{z+1}{2z+1})-\phi(\frac{z}{2z+1})) \]
and hence
\[ \phi(-z-1)-\phi(-z)=-(\phi(z)-\phi(z+1)=\phi(z+1)-\phi(z)\]
Again decomposing $\phi$ in an even and an odd part we see that this
implies $\phi\in\P^{e}_{n+1}$.  Now note that
\[  \tilde T^{-1}\Upsilon_{n}=\{\phi\in\P_{n+1}|\tilde T\phi(z)\in \Upsilon_{n}\}\]
\[ =\{\phi\in\P_{n+1}|\phi(z)-\phi(z+1)=-(-2z-1)^{n}(\phi(-\frac{z+1}{2z+1})-\phi(\frac{z}{2z+1}))\}    \] 
\[  =\{\phi\in\P_{n+1}|\phi(z)-\phi(z+1)=(2z+1)^{n}(-\phi(-\frac{z+1}{2z+1})+\phi(\frac{z}{2z+1}))\}  \]
>From this equation and the fact that polynomials in
$\tilde\wp^{+}_{n}$ are odd and polynomials in $\tilde\wp^{-}_{n}$ are
even we get
\[ \tilde\wp_{n}=\tilde \wp^{+}_{n}\oplus \tilde\wp^{-}_{n}\subseteq \tilde T^{-1}\Upsilon_{n}\]
It remains to show that
\[  \tilde T^{-1} \Upsilon_{n}\subseteq \tilde \wp^{+}_{n}\oplus \tilde\wp^{-}_{n}\]
Assume that there exists $\phi\in \tilde T^{-1} \Upsilon_{n}\backslash
(\tilde \wp^{+}_{n}\oplus \tilde\wp^{-}_{n})$. Decompose
$\phi=\phi^{e}+\phi^{o}$ with $\phi^{e}\in\P_{n+1}^{e}$ and
$\phi^{o}\in\P_{n+1}^{o}$.  Define $\tilde \phi$ by
\[ \tilde\phi(z):=\phi^{e}(z)-\phi^{e}(z+1)-(2z+1)^{n}(-\phi^{e}(\frac{z+1}{2z+1})+\phi^{e}(\frac{z}{2z+1})).\]
Since we assumed $\phi\in \tilde T^{-1} \Upsilon_{n}$ we have
\[ -\tilde\phi(z)=\phi^{o}(z)-\phi^{o}(z+1)-(2z+1)^{n}(\phi^{o}(\frac{z+1}{2z+1})+\phi^{o}(\frac{z}{2z+1})) \]
Substituting $-z-1$ for $z$ in these equations yields
\[\tilde\phi(z)=-\tilde\phi(-z-1)\qquad\mbox{and}\qquad
\tilde\phi(z)=\tilde\phi(-z-1).\] But this obviously implies
$\tilde\phi=0$.  On the other hand since $\phi\not\in (\tilde
\wp^{+}_{n}\oplus \tilde\wp^{-}_{n})$ the polynomial $\tilde\phi$ can
not be zero.  This is a contradiction and our proof is complete.
\qed\\
Since the kernel of $\tilde T$ consists of all constant polynomials
which are contained as well in $\tilde\wp_{n}^{-}$ Proposition 8.1 and
section 3 has the following corollary.
\begin{corollary}
  For all $n\ge 2$ even we have
\[ \dim\tilde\wp_{n}=\dim\Upsilon_{n}+1=2[(k-1)/4]+3.\]
\end{corollary}
Combining the following Proposition with this Corollary we see that
Conjecture 8.1 is really true in the case $N=2$.
\begin{proposition}
  For all $n\ge 2$ even we have
\[ \dim\tilde\wp_{n}=\dim \wp_{n}(\Gamma_{0}(2))+1 \] 
\end{proposition}
\proof Looking again at the proof of Proposition 8.1 we see that
\[ \wp_{n}(\Gamma_{0}(2))\cong \{\phi\in\P_{n}|\phi(z)-\phi(z+1)=(2z+1)^{n}(-\phi(-\frac{z+1}{2z+1})+\phi(\frac{z}{2z+1}))\}\]
and we will identify these isomorphic spaces in the following.\\
We can write every $\tilde\phi\in\tilde\wp_{n}$ as
$\tilde\phi=\phi+cz^{n+1}$ where $\phi\in\P_{n}$ and $c$ is a constant
such that
\[ \phi(z)-\phi(z+1)+(2z+1)^{n}(\phi(-\frac{z+1}{2z+1})-\phi(\frac{z}{2z+1}))\]\[ =c(z^{n+1}-(z+1)^{n+1}-\frac{(z+1)^{n+1}}{2z+1}-\frac{z^{n+1}}{2z+1}) \]
Let
\[ \check\phi(z)=z^{n+1}-(z+1)^{n+1}-\frac{(z+1)^{n+1}}{2z+1}-\frac{z^{n+1}}{2z+1} \]
Since $(z^{n+1}+(z+1)^{n+1})$ has a zero at $z=-1/2$ if $n$ is even we
have $\check \phi\in\P_{n}$.  Now there exists a $\phi_{1}\in
\P_{n}$ such that
\[ \phi_{1}(z)-\phi_{1}(z+1)+(2z+1)^{n}(\phi_{1}(-\frac{z+1}{2z+1})-\phi_{1}(\frac{z}{2z+1}))=\check \phi(z)\]\
Hence we have
\[ \tilde\wp_{n}\cong\{\phi+cz^{n+1}|\phi\in\P_{n},\quad c\in\R,\quad  \phi
\in\wp_{n}(\Gamma_{0}(2))+c\phi_{1}\}\] \[\cong
\wp_{n}(\Gamma_{0}(2))\oplus <\phi_{1}+z^{n+1}> \] which proves our
Proposition.
\qed\\\\
{\bf The case $N>2$}\\\\
In the case $N>2$ prime we have some numerical evidence for our
general conjecture.  Table 2 below contains the dimension of
$\wp_{n}(\Gamma_{0}(N))$ for some values $N>2$ prime and $n\ge 2$
even.  If we compare Table 2 with $\dim E_{k}(\Gamma_{0}(N))$ using
the table in the appendix of $\cite{[MI]}$ for $\dim
C_{n}(\Gamma_{0}(N)$ and Proposition 3.2 we see that Conjecture 8.1 is
true for all values
of $N$ and $n$ we have checked.\\
\begin{center}
\begin{tabular}{|c|c|c|c|c|c|c|}
\hline
      ~&  $N=3$ & $N=5$ & $N=7$ & $N=11$ & $N=13$ & $N=17$    \\      
\hline
n=2 & 2 &4 & 4 & 6 & 8 & 10    \\ 
\hline
n=4 & 4 & 4 & 8 & 10 & 12           & 14   \\ 
\hline
n=6 &  4 & 8  & 8 & 14 &  16           & 22  \\ 
\hline
n=8 &  6 & 8 & 12& 18 &  20      & 26  \\ 
\hline
n=10 &  8 & 12  &16  &  22   & 28 & 34    \\ 
\hline
n=12  &  8  & 12  & 16 & 26   & 28   & 38    \\ 
\hline
n=14  &  10  & 16  & 20 & 30   & 36   & 46   \\ 
\hline
n=16  &  12  & 16  & 24 & 34   & 40   & 50   \\ 
\hline
n=18  &  12  & 20  & 24 & 38   & 44   & 58  \\ 
\hline
n=20  &  14 & 20  & 28 & 42   & 48   & 62  \\ 
\hline
n=22  &  16 & 24  & 32 & 46  & 56 & 70  \\ 
\hline
n=24  &  16  & 24  & 32 & 50   & 56  & 74   \\ 
\hline
\end{tabular}
\end{center}
{\bf Table 2: $\dim\wp_{n}(\Gamma_{0}(N))$}\\\\
We have calculated these dimensions algebraically using Matematica.
If $N$ is prime we have the simple system of representatives of the
coset set $\tilde\Gamma_{0}(N)$ given in section 2.  This allows us to
calculate the action of $Q$ and $T$ on the coset set and thus the
induced representation $\chi^{\Gamma}$. Now we proceed in exactly the
same way as in the calculations for $\Gamma(N)$.  Here we have to
solve a system of only $(N+1)n$ linear equations.


\begin{thebibliography}\small
  \bibitem{[CM1]} C.-H. Chang and D. Mayer, {\it The transfer
    operator approach to Selberg's zeta function and modular and Maass
    wave forms for $PSL(2,\Z)$}, in D. Hejhal and M. Gutzwiller et al,
  editors, IMA Volumes {\bf 109} 'Emerging applications of number theory',
  Springer Verlag 72-143 (1999).  \bibitem{[CM2]} C.-H. Chang and
  D. Mayer, {\it Thermodynamic formalism and Selberg's zeta function
    for modular groups}, Regular and Chaotic Dynamics {\bf 5}, 281-312
  (2000).  \bibitem{[CM3]} C.-H. Chang and D. Mayer, {\it
    Eigenfunctions of the transfer operators and the period functions
    for modular groups}, in 'Proceedings of an AMS workshop on
  arithmetic, spectral and dynamical zeta functions', San Antionio
  (1999); to appear in: Contemporary. Math (2001).  \bibitem{[EI]} M.
  Eichler, {\it Eine Verallgemeinerung der abelschen Integrale} ,
  Math. Zeitschrift,  {\bf 67}, 267-298 (1957).  \bibitem{[FH]} W.
  Fulton and J. Harris, {\it Representation Theory - A first Cours},
  Springer Verlag (1991).  \bibitem{[KH]} A. Katok and B. Hasselblatt,
  {\it Modern Theory of Dynamical Systems}, Cambridge University press
  (1995).  \bibitem{[LM]} A. Lasota and M. Mackey, {\it Probabilistic
    properties of deterministic systems}, Cambridge University Press,
  (1985).  \bibitem{[LW]} J. Lewis, {\it Spaces of holomorphic
    functions equivalent to Maass cusp forms}, Invent. Math., {\bf 127(2)},
  271-306 (1997).  \bibitem{[LZ]} J. Lewis and D. Zagier, {\it
    Period functions and the Selberg's zeta function for modular
    groups}, in The Mathematical Beauty of Physics, Adv. Series in
  Math. Physics 24, World Scientific, Singapour, 83-97  (1997).
\bibitem{[MA]} D. Mayer, {\it The thermodynamic formalism approach to
    Selberg's zeta function for $PSL(2,\Z)$}, Bull. Am. Math. Soc.,
  {\bf 25}, 55-60 (1991).  \bibitem{[MI]} T. Miyake, {\it Modular Forms},
  Springer Verlag (1989).  \bibitem{[SCH]} G. Shimura, {\it
    Introduction to the Arithmetic Theory of Automorphic Functions},
  Princeton University Press (1971).  \bibitem{[VE]} A.B. Venkov, {\it
    Spectral Theory of Automorphic Functions and Its Applications},
  Kluwer Academic Publishers (1990).
\end{thebibliography}
\end{document}